%% file: manuscript.tex
\documentclass[reqno]{amsart}

\usepackage[margin=3cm]{geometry}

\usepackage{arxiv_style}
\usepackage{fontawesome}

\makeatletter
\def\blfootnote{\gdef\@thefnmark{}\@footnotetext}
\makeatother

\newcommand{\fakesubsection}[1]{\medskip\noindent\textbf{#1.}}

\title[Structural Identifiability and Comparative Calibration of WRCs in Porous Media]{Structural Identifiability and Comparative Calibration of Water Retention Curves for Imbibition in Porous Media}

\author[G.\ Bretti, M.\ Ceseri, E.\ Onofri, M.\ Paoluzzi]{}

\begin{document}

\blfootnote{$^{\star}$ Corresponding Author, (\href{mailto:elia.onofri@kaust.edu.sa}{\faEnvelopeO}) \texttt{elia.onofri@kaust.edu.sa}, (\href{https://www.eliaonofri.it}{\faGlobe}) \texttt{www.eliaonofri.it}}

\maketitle

\vspace{-1em}

\begin{center}
    \begin{minipage}{.89\linewidth}\centering
        \textsc{Gabriella Bretti}$^{\,\textsc{a},\, \orcidlink{0000-0001-5293-2115}}$,
        \textsc{Maurizio Ceseri}$^{\,\textsc{a},\, \orcidlink{0000-0002-7055-9323}}$,
        \textsc{Elia Onofri}$^{\,\textsc{b},\, \textsc{a},\, \star,\, \orcidlink{0000-0001-8391-2563}}$, \\
        \textsc{Matteo Paoluzzi}$^{\,\textsc{c},\, \orcidlink{0000-0003-3983-8161}}$.
        \\
        \bigskip
        \begin{minipage}{.45\linewidth}\centering
            \footnotesize
            $^\textsc{a}$Istituto per le Applicazioni del Calcolo (IAC),\\
            Consiglio Nazionale delle Ricerche (CNR)\\
            Rome 00185, Italy
        \end{minipage}
        \begin{minipage}{.45\linewidth}\centering
            \footnotesize
            $^\textsc{c}$Department of Physics,\\
            Sapienza, Università di Roma\\
            Rome 00185, Italy
        \end{minipage}
        
        \bigskip

        \begin{minipage}{.9\linewidth}\centering
            \footnotesize
            $^\textsc{b}$Computer, Electrical and Mathematical Sciences and Engineering (CEMSE) Division,\\
            King Abdullah University of Science and Technology (KAUST)\\
            Thuwal 23955, Saudi Arabia
        \end{minipage}
    \end{minipage}
\end{center}

\medskip
\thispagestyle{empty}

\begin{abstract}

    This paper investigates the structural identifiability and a comparative calibration of four water retention curves (WRCs) within the framework of the Richards equation coupled with Darcy's law for capillary imbibition in porous media. The considered models --two classical physically-based laws and two abstract parametrisations developed for building stones-- are consistently reformulated by expressing the hydraulic conductivity $K(\Theta)$ and capillary pressure $\psi(\Theta)$ independently, allowing the nonlinear diffusion coefficient $D(\Theta)$ to be reconstructed in a unified structural form. This common representation enables a rigorous mathematical comparison across models with different theoretical foundations.

    \noindent All models are calibrated against the same experimental imbibition dataset using a grid-based optimisation strategy with adaptive refinement. The analysis reveals a structural property of the associated inverse problem: the hydraulic conductivity and the capillary pressure scale enter the governing equation multiplicatively and therefore cannot be independently identified from imbibition data. Only their product acts as an observable diffusion parameter, where model discrimination is primarily governed by the shape of the resulting effective diffusion function.

    \noindent To the best of our knowledge, this is the first study providing a coherent cross-calibration of these WRCs against an identical dataset within a unified computational framework. Our open-source implementation, released within the Stoneverse platform, provides a reproducible baseline for further developments, including probabilistic inversion and learning-based approaches.
    
    \medskip
    
    
    \smallskip
    
    \noindent{\bf AMS-MSC 2020:} 35K65 \sep 76S05.

\end{abstract}

\begin{multicols}{2}
    \tableofcontents
\end{multicols}

\newpage


\section{Introduction}

Water flow in partially saturated porous media is commonly described by the Richards equation, a nonlinear parabolic equation coupling mass conservation with Darcy’s law. This model plays a central role in hydrology, soil science, and material science, and is particularly relevant for the study of moisture transport in natural and artificial stones.

\begin{figure}[h]
    \centering
    \includegraphics[scale=0.25]{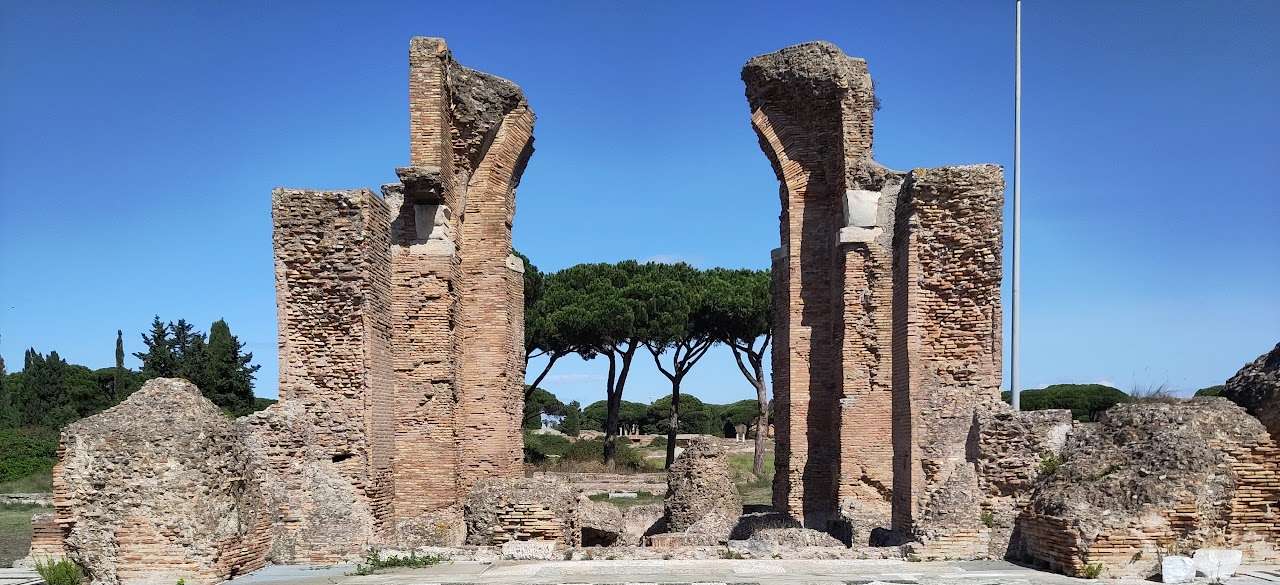}
    \caption{Roman pillars made of bricks in the archaelogical park of Ostia Antica.}
    \label{fig:ostia}
\end{figure}

In architectural heritage materials, capillary imbibition represents one of the primary mechanisms driving deterioration. Porous stones exposed to environmental agents (see Figure~\ref{fig:ostia}) undergo degradation processes such as salt crystallisation, freeze–thaw cycles, and dissolution phenomena, all strongly influenced by moisture ingress. Monitoring and modelling capillary absorption are therefore essential for understanding long-term material behaviour.

In one spatial dimension and neglecting gravity --an assumption justified for small laboratory specimens-- the Richards equation reads
\begin{equation}\label{eq:Richards-pres}
    \frac{\partial\theta(\psi)}{\partial t}=\frac{\partial}{\partial z}\left(K(\psi)\frac{\partial\psi}{\partial z}\right),
\end{equation}
where $\theta$ denotes the volumetric water content, $K$ the hydraulic conductivity, and $\psi$ the pressure head. The pressure head is negative in the unsaturated regime and positive under saturated conditions. The hydraulic conductivity is a positive, increasing function of $\psi$ in the unsaturated range and becomes constant at saturation. The water content satisfies $\theta \in [\theta_r,\theta_s]$, where $\theta_r$ is the residual water content and $\theta_s$ the saturated water content.
Do note that gravity is neglected because the laboratory specimens considered in the experimental setting are sufficiently small for gravitational effects to be negligible compared with capillary forces.

To close equation \eqref{eq:Richards-pres}, a constitutive relation between $\theta$ and $\psi$ must be specified.
This relation is commonly referred to as the water retention curve (WRC).
The choice of the WRC directly influences the nonlinear behaviour of the model.

The Richards equation can be reformulated in diffusion form as
\begin{equation}\label{eq:Richards}
\frac{\partial\theta}{\partial t}=\frac{\partial}{\partial z}\left(D(\theta)\frac{\partial\theta}{\partial z}\right)
\end{equation}
where $D(\theta)$ is the nonlinear diffusion function defined by:
\begin{equation}\label{eq:diffusion}
D(\theta)=K(\psi(\theta))\frac{\partial \psi}{\partial \theta}.
\end{equation}
In this representation, the WRC determines the structure of the nonlinear diffusion coefficient $D(\theta)$. Different choices of retention law therefore lead to diffusion operators with substantially different degeneracy, curvature, and scaling properties.

A wide variety of parametric WRCs has been proposed in the literature. In this work we consider four representative models: the Brooks–Corey (BC, \cite{brooks-corey}) and van Genuchten (vG, \cite{van-genuchten}) laws, widely used in soil physics, and two formulations developed in the context of building stones, namely the Natalini–Nitsch (NN, \cite{bracciale,goid,Clarelli2010}) model and the more general Bretti–et al.\ (BkP~\cite{malte-catania}) family.

The aim of this study is twofold. First, we provide a unified mathematical framework in which all four models are consistently formulated by independently specifying $K(\Theta)$ and $\psi(\Theta)$ and reconstructing the associated diffusion function $D(\Theta)$. Second, we perform a systematic comparative calibration of these models against an identical imbibition dataset, thereby analysing both their predictive behaviour and the structural properties of the corresponding inverse problem.


\subsubsection*{Main contributions}

The main contributions of this work can be summarised as follows:

\begin{itemize}
    \item[(i)] A unified modelling framework for four water retention curves is introduced by consistently reconstructing the nonlinear diffusion coefficient from independently specified conductivity and capillary pressure functions within the Richards equation under Darcy's Law.
    
    \item[(ii)] A systematic cross-calibration of the models against an identical imbibition dataset is performed, enabling a direct structural comparison under homogeneous numerical and optimisation settings.
    
    \item[(iii)] A structural identifiability analysis of the associated inverse problem is provided, showing that only the product of the hydraulic conductivity scale and the capillary pressure scale is observable in imbibition experiments.
\end{itemize}


\subsubsection*{Roadmap}

The rest of the paper is organised as follows.
Section \ref{sec:literature} provides a short overview over the literature context.
Section \ref{sec:wrc} presents the overall framework and the four water retention models.
Section \ref{sec:exp-setting} describes the experimental setting, the numerical discretisation of the direct problem and the calibration procedure of the inverse problem.
Section~\ref{sec:res} presents the comparisons results.
Finally, Section \ref{sec:concl} concludes the paper discussing the finding, highlighting the limitations and outlining a few future research directions.


\section{Literature Context}\label{sec:literature}

The mathematical modelling of liquid transport in variably saturated porous media has been extensively investigated over the past decades, primarily within the framework of the Richards equation. 
A central modelling component is the specification of constitutive relations linking volumetric water content, capillary pressure, and hydraulic conductivity. 
These relations determine the nonlinear structure of the governing parabolic equation and strongly influence both qualitative solution behaviour and numerical complexity.
Consequently, extensive research has been devoted to modelling water and solute transport in variably saturated soils and stones \cite{ipp, scelsi, scarfone, zahasky, suh}, covering both experimental characterisation and computational approaches.  

Classical water retention curves (WRCs) were developed by introducing parametrisations connected to pore size distribution. 
Although phenomenological at the macroscopic scale, such models implicitly encode microstructural effects through a small number of parameters \cite{brooks-corey, van-genuchten, fredlund}. 
In this setting, the volumetric water content $\theta$ is expressed as a function of the pressure head $\psi$, and the unsaturated hydraulic conductivity $K(\theta)$ is derived accordingly.
These parametric relations are typically used either to interpolate experimental data or to estimate hydraulic properties via inverse modelling \cite{celia1990}.

A significant contribution to this framework is due to Mualem \cite{mualem}, who proposed a theoretical model for the relative hydraulic conductivity associated with a given retention curve. 
Since direct measurement of unsaturated conductivity is experimentally challenging, Mualem's formulation has been widely combined with different WRC parametrisations, most notably the Brooks--Corey (BC) \cite{brooks-corey} and van Genuchten (vG) \cite{van-genuchten} models.

The BC model introduces an explicit air-entry pressure corresponding to the largest pore size, resulting in a piecewise-defined retention behaviour with a sharp transition to saturation. 
In contrast, the vG model provides a smooth and continuous relation without an explicit air-entry threshold. 
Variants incorporating an air-entry correction have also been proposed, for instance by Vogel et al.~\cite{vogel}, in order to better describe behaviour near saturation.

Beyond soil-oriented parametrisations, alternative formulations have been introduced in the context of moisture transport in building stones. 
In \cite{Clarelli2010}, the Natalini--Nitsch (NN) model was proposed as a structurally unified description in which hydraulic conductivity and capillary pressure are both expressed as functions of saturation. 
In this framework, flow is restricted to a closed saturation interval: permeability vanishes below residual saturation, while capillary pressure approaches zero at maximal saturation.

More recently, Bretti et al.~\cite{malte-catania} proposed a broader family of absorption functions (BkP) in which the contributions of capillary pressure and conductivity are separated and additional parameters control curvature and degeneracy properties. 
This increased flexibility allows a wider range of diffusion behaviours within the Richards framework.

While these models have been widely studied individually, comparative analyses are typically performed under heterogeneous modelling assumptions or distinct numerical settings. 
Moreover, the structural implications of these parametrisations on the associated inverse problem have received limited attention. 
In particular, the identifiability of scaling parameters arising from the multiplicative interaction between conductivity and capillary pressure is rarely addressed explicitly.

The present work contributes to this gap by reformulating four representative WRCs --BC, vG, NN, and BkP-- within a unified diffusion framework, enabling a consistent cross-calibration against an identical dataset and a structural comparison of their inverse properties.


\section{Water retention curves}\label{sec:wrc}

Before introducing any specific WRC, it is convenient to reformulate equation \eqref{eq:Richards} in terms of the effective saturation $\Theta \in [0, 1]$
\begin{equation}\label{eq:sat_def}
    \Theta = \frac{\theta-\theta_r}{\theta_s-\theta_r},\mbox{ i.e. } \theta=(\theta_s-\theta_r)\Theta+\theta_r.
\end{equation}
With this new variable, Richards equation becomes:
\begin{equation}\label{eq:Richards_sat}
\frac{\partial\Theta}{\partial t}=\frac{\partial}{\partial z}\left(D(\Theta)\frac{\partial\Theta}{\partial z}\right)
\end{equation}
with the non-linear diffusion function given by the Darcy' law that reads
\begin{equation}\label{eq:diffusion_sat_new}
D(\Theta)=\frac{K(\psi(\Theta))}{\mu(\theta_s-\theta_r)}\frac{\partial \psi}{\partial \Theta}\,,
\end{equation}
where $\mu$ represent the liquid viscosity ($\mu_\textsc{water}=8.9$ millipoise at 25 $^\circ$C).
Equation \eqref{eq:Richards_sat} is coupled with null initial condition and boundary conditions of Robin-type:
\begin{eqnarray}\label{DBC}
\Theta(x,t)&=&0, \textrm{for} \ x>0 \nonumber\\
\Theta(0,t)&=&1,\textrm{for} \ t\ge 0,\\
\frac{\partial\Theta(L,t)}{\partial z} &=&K_w (\Theta_\textrm{ext}-\Theta(L,t)) ,\textrm{for} \ t\ge 0 \nonumber.
\end{eqnarray}
These conditions model a specimen that is initially dry, fully saturated at the lower boundary due to contact with a water source (sponge), and exchanging moisture at the upper boundary with the external environment at rate $K_w$. Following \cite{Stolfi25}, we set $K_w=10^{\klog}$, with $\klog \in \mathbb{R^-}$ and treat $\klog$ as a calibration parameter.


\subsection{General framework for WRCs}

We define a water retention curve as the constitutive relation between pressure head and saturation,
\begin{equation}\label{eg:general_wrc}
    \Theta=\Theta(\psi)
\end{equation}

In the unsaturated regime ($\psi \le 0$), $\Theta(\psi)$ is monotone increasing, while at saturation ($\psi>0$) one has $\Theta = 1$.
Under suitable regularity assumptions, the relation can be inverted for $\Theta \neq 1$, yielding $\psi = \psi(\Theta)$.

Once the retention curve is specified, the hydraulic conductivity can be derived. A widely used approach is due to Mualem \cite{mualem}, who first defined the hydraulic conductivity as
\begin{equation}\label{eq:mualem}
    K(\Theta)=K_s\Theta^{1/2}\left[\frac{\int_0^\Theta \frac{ds}{\psi(s)}}{\int_0^1 \frac{ds}{\psi(s)}}\right]^2
\end{equation}
where $K_s$ denotes the hydraulic conductivity at saturation.

The diffusion coefficient $D(\Theta)$ is then obtained through \eqref{eq:diffusion_sat_new}.
Different WRCs therefore induce different nonlinear diffusion structures.

In the following, we describe the four models considered in this work.


\subsection{Brooks-Corey (BC)}

One of the most applied WRC is due to Brooks and Corey~\cite{brooks-corey}.
This formulation assumes a power-law relation between pressure head and effective saturation and is given by:
\begin{equation}\label{BC}
    \Theta=\left(\frac{\psi}{\psi_b}\right)^{-\lambda},\mbox{ or }\psi=\psi_b\Theta^{-1/\lambda}
\end{equation}
where $\psi_b < 0$ is the bubbling (air-entry) pressure head and $\lambda>0$ is a pore-size distribution parameter (nondimensional).
By applying Mualem's theory, we can determine the hydraulic conductivity as:
\begin{equation}
    K(\Theta)=K_s\Theta^{\frac{2}{\lambda}+\frac{5}{2}}
\end{equation}
and from equation \eqref{eq:diffusion_sat_new}
\begin{equation}\label{BC_diff}
    D(\Theta)=D_0\Theta^{\frac{1}{\lambda}+\frac{3}{2}}\,,
\end{equation}
with
\begin{equation}\label{D0}
    D_0=-\frac{K_s\psi_b}{\mu(\theta_s-\theta_r)} \frac{1}{\lambda}.
\end{equation}
Thus, the BC model leads to a pure power-law diffusion function.


\subsection{van Genuchten (vG)}

The van Genuchten \cite{van-genuchten} retention curve is defined by
\begin{equation}\label{MvG}
    \Theta=\left[1+\left(\frac{\psi}{\psi_b}\right)^n\right]^{-\frac{\lambda}{n}}\,,\quad\text{or}\quad\psi=\psi_b(\Theta^{-\frac{n}{\lambda}}-1)^{1/n}\,,
\end{equation}
with $n$ being another fitting parameter, while $\lambda$ (and $\psi_b$) defined as in the BC curve\footnote{It is important to note that the original notation for \eqref{MvG} from \cite{van-genuchten} adopted $\alpha = \sfrac{1}{\psi_b}$ instead and $m$ instead of $\sfrac\lambda n$. In this version, we decided to keep concistency with BC formulation.}.
Using Mualem's theory and imposing the classical contraint $n=\lambda+1$ we obtain a closed-form expression of the hydraulic conductivity
\begin{equation}
    K(\Theta)=K_s\Theta^{1/2}[1-(1-\Theta^{\sfrac1m})^{m}]\,, \qquad m = \frac{\lambda}{n} = \frac{n-1}n = \frac{\lambda}{\lambda+1}.
\end{equation}

Plugging \eqref{MvG} into \eqref{eq:diffusion_sat_new} leads to the following expression of the nonlinear diffusion function:
\begin{equation}\label{MvG_diff}
    D(\Theta)=D_0\frac{\left[1-(1-\Theta^{\sfrac1m})^{m} \right]^2(1-\Theta^{\sfrac1m})^{-m}}{\Theta^{1/\lambda+1/2}}
\end{equation}
with $D_0$ as per \eqref{D0}, or, analogously
\begin{equation}
    D(\Theta)=D_0 \frac{(1-\omega)^2}{\omega} \Theta^{\sfrac12 - \sfrac1m}, \qquad \omega = (1-\Theta^{\sfrac1m})^{m} \ .
\end{equation}

\begin{remark}
    The constraint $n=\lambda+1$ ensures closed-form integrals in \eqref{eq:mualem}, thereby reducing the number of independent calibration parameters. If $n$ and $\lambda$ are treated as independent, the integrals must be evaluated numerically and both parameters must be calibrated separately.
\end{remark}
\begin{remark}
    For $|\psi| \to \infty$ (low saturation), the van Genuchten model asymptotically reduces to a Brooks–Corey-type power law.
\end{remark}

\subsection{Natalini--Nitsch (NN)}
In the NN model, originally introduced in \cite{Clarelli2010}, the diffusion function is prescribed directly as
\begin{equation}\label{NN_diff}
    D(\Theta)=D_0\Theta(1-\Theta)
\end{equation}
with 
\begin{equation}
D_0=-\frac{K_s\psi_b}{\mu(\theta_s-\theta_r)}.
\end{equation}

In this case, the associated constitutive functions can be chosen as
\begin{equation}
K(\Theta)=K_s\Theta^4
\end{equation}
and we can assume the following form of the capillary pressure head
\begin{equation}\label{NN}
    \psi(\Theta)=\frac{\psi_b}{2}\frac{(1-\Theta)^2}{\Theta^2}.
\end{equation}
This formulation separates permeability and capillary pressure in a way that was not explicitly adopted in the original presentation of \cite{Clarelli2010}, where the diffusion function was introduced directly.


\subsection{BkP model}

The BkP model \cite{malte-catania} generalises the NN formulation by explicitly separating permeability and capillary contributions since its conception.
This allows derive specific parameters related to the capillary pressure and to the permeability properties of the porous material that can be obtained experimentally. 
The hydraulic conductivity is assumed to be
\begin{equation}\label{K_BCP}
    K(\Theta)=K_s\Theta^\gamma
\end{equation}
with $K_{s}>0$ being the conductivity at saturation as usual and $\gamma>0$ a fitting parameter.

The capillary pressure is defined as:
\begin{equation}\label{BCP}
    \psi (\Theta) = \psi_b(\theta_s-\theta_r)^{2-\alpha}\frac{(1-\Theta)^2}{\Theta^{\alpha}}
\end{equation}
with $\psi_b$ the bubbling pressure and $\alpha>0$ a fitting parameter.

Then, the resulting diffusion function yields a family of curves given by:
\begin{equation}\label{BkP_diff}
    D(\Theta)=D_0(1-\Theta)\Theta^{\gamma-\alpha}\left(2+\alpha\frac{1-\Theta}{\Theta}\right),
\end{equation}
where
\begin{equation}
    D_0=-\frac{K_s\psi_b}{\mu}(\theta_s-\theta_r)^{1-\alpha}.
\end{equation}

\begin{remark}
    It is worth noting that \eqref{BkP_diff} reduces to \eqref{NN_diff} when $\alpha=2$ in \eqref{BCP} and $\gamma=4$.
\end{remark}


\subsection{Parameter sets}\label{sec:parameter_set}

For the numerical analysis we later present in Section~\ref{sec:res}, we adopt the following parameterisation strategy.
\begin{description}
    \item[Porosity-related parameters]
    The residual and saturated water contents are expressed as fractions of the open porosity $n_0$:
    \begin{equation}\label{eq:a_b}
        a = \frac{\theta_r}{n_0},
        \qquad
        b = \frac{\theta_s}{n_0}.
    \end{equation}
    The open porosity is assumed known from laboratory measurements: in what follows, we adopt $n_0=0.285$, corresponding to a common brick \cite{bracciale}.

    \item[Boundary exchange parameter] The parameter $\klog=\log(K_w)$ is calibrated in order to assess the impact of boundary exchange on imbibition dynamics.

    \item[Model-specific shape parameters]
    These include $\lambda$ for BC and vG, the couple $(\alpha,\gamma)$ for BkP, and no additional parameters for NN, since $D(\Theta)$ is fixed up to $D_0$.

    \item[Maximum diffusion scaling]
    To enable a unified parametrisation across models with structurally different retention functions, we introduce the combined scaling parameter
    \begin{equation}\label{eq:c}
        c = K_s \, \psi_b \,,
    \end{equation}
    which represents the characteristic diffusion strength of the material;
    Here, it is important to note that in many formulations, the viscosity term $\mu$ is also absorbed into $c$.
    Indeed, in the Richards-type formulation adopted here, cf.~\eqref{eq:Richards}, the effective diffusion coefficient can be written in the generic form
    \begin{equation}\label{eq:effective_D}
        D(\Theta) = c \cdot \widehat{D}(\Theta),
    \end{equation}
    where $\widehat{D}(\Theta)$ is a dimensionless function determined by the specific retention law.
    As a consequence, the parameter $c$ governs the overall magnitude of moisture transport, while the shape of $D(\Theta)$ is controlled by the model-specific parameters.
\end{description}

\begin{remark}
    It is worth noticing that, since $K_s$ and $\psi_b$ appear multiplicatively in the diffusion term, their individual contributions cannot be independently identified from imbibition data alone, and only their product is structurally observable in this experimental setting.
    In what follow, for the sake of comparability, we set $K_s = 1$ and we derive $\psi_b$ as a function of $K_s$ for each model, namely, $\psi_b = \sfrac{c}{K_s} = c$.
\end{remark}


\section{Experimental benchmark and forward simulation}\label{sec:exp-setting}

This section presents the assessment for the predictive performance of the different water retention formulations by comparison with a controlled imbibition experiments.
The objective is here to determine whether distinct constitutive laws, once calibrated to the same dataset, produce comparable macroscopic responses and internal moisture dynamics and whether common parameters (like saturated and residual water content) fall in comparable intervals.

To this end, we first introduce the experimental benchmark and define the measured observable.
We then describe the forward problem and its numerical approximation, which together define the mapping from a parameter vector to the simulated imbibition curve.
The parameter identification procedure is presented in the subsequent subsection.


\subsection{Experimental benchmark}

\begin{figure}[t]
    \centering
    \begin{subfigure}{.35\linewidth}
        \centering
        \scalebox{.9}{\input{figs/experimentalSetting}}
        \caption{}
        \label{fig:benchmark-a}
    \end{subfigure}
    \begin{subfigure}{.6\linewidth}\centering
        \includegraphics[width=.9\linewidth, clip, trim={0cm 0.55cm 0cm 0cm}]{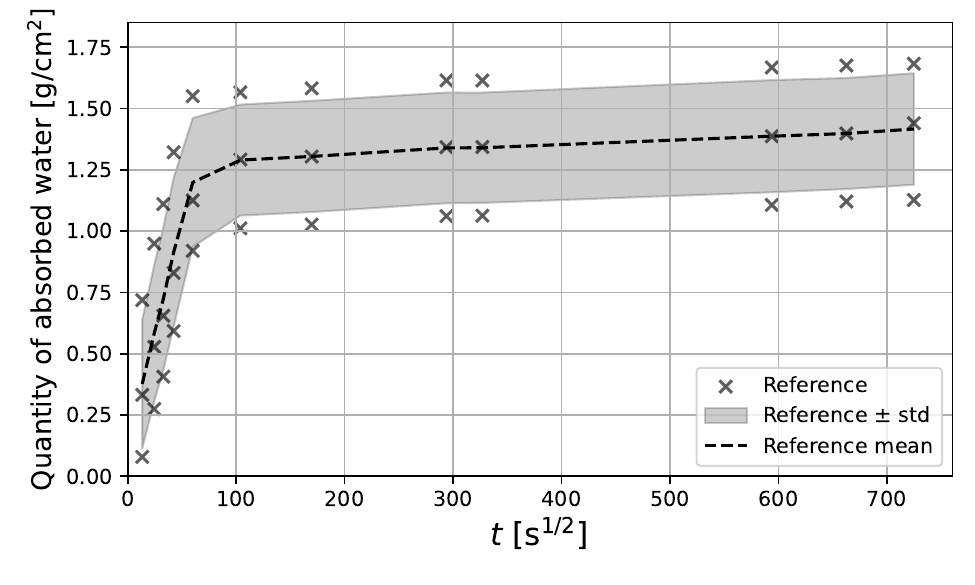}
        \caption{}
        \label{fig:benchmark-b}
    \end{subfigure}
    \caption{
        Experiment from~\cite{bracciale} used to to calibrate the models. 
        \subref{fig:benchmark-a} Depiction of the experimental setup. 
        \subref{fig:benchmark-b} Imbibition curves obtained from the common brick.
        The three series of reference data (grey crosses) are averaged to create a single reference curve (dashed line).
    }
    \label{fig:benchmark}
\end{figure}

The experimental dataset (taken from \cite{bracciale}) consists of one-dimensional capillary imbibition tests (capillary rise) performed on cubic brick samples with edge length $L = 5\,\mathrm{cm}$.
One face of each specimen was placed in contact with water in a closed environment, while the lateral faces were sealed to enforce approximately one-dimensional moisture transport.
Top face was left open to allow vapour exchange with the closed environment, which was forced into high water saturation levels ($\theta_\textrm{ext}= 0.06254$\,g/cm$^3$, as in \cite{bracciale}).
Figure~\ref{fig:benchmark-a} present a visual depiction of the experiment.

Water uptake was measured gravimetrically at prescribed time intervals $t_0, \dots,$ $t_{N_T}$.
The absorbed mass was converted into water uptake per unit surface area, yielding the experimental imbibition curve $Q(t)$, sampled as a vector of measurements $\mathbf{Q} = (Q_0, \dots, Q_{N_T})$, where
\begin{equation}
    Q(t_k) = Q_k, \quad k = 0,\dots,N_T \,.
\end{equation}
The averaged curve, depicted in Figure~\ref{fig:benchmark-b} is used as reference data in the calibration procedure.
\begin{remark}  
    To reduce experimental variability, three tests were performed under identical conditions and the resulting mass readings were averaged.
\end{remark}

For the numerical simulations, the physical domain is represented as a one-dimensional interval of length $L$ corresponding to the sample thickness.
We hence set $z \in [0,L]$, where $z=0$ represents the water-exposed face and $z=L$ the open boundary.

The experimentally measured quantity to be reproduced by the models is therefore the cumulative water uptake per unit surface area as a function of time.


\subsection{Forward problem and numerical approximation}


\subsubsection*{Governing equation}

Moisture transport in the porous medium is described in terms of the volumetric water content $\Theta(z,t) \in [0, 1]$, where we recall from our assumption \eqref{eq:sat_def} that $\Theta(z, t) = \frac{\theta(z, t) - \theta_r}{\theta_s - \theta_r}$.
In one spatial dimension, the governing equation reads
\begin{equation}\label{eq:diff-eq}
	\partial_t \Theta = \partial_z \bigl( D(\Theta) \, \partial_z \Theta \bigr),
\end{equation}
where $D(\Theta)$ is the moisture diffusivity, determined by the chosen WRC and hydraulic conductivity model.

The boundary and initial conditions correspond to the imbibition setup \eqref{DBC}:
\begin{equation}\label{IBC}
	\begin{cases}
		\Theta(z,0) = \Theta_0	& \forall z \in (0, L]\\
		\Theta(0,t) = 1			& \forall t \in [0, t_{N_T}]\\
		\partial_z \Theta(L,t) = K_w \left( \Theta_\text{ext} - \Theta(L,t) \right)  & \forall t \in [0, t_{N_T}]\\
	\end{cases}
\end{equation}
or, in other words, that (i) the specimen is assumed initially at uniform moisture content $\Theta_0 \sim 0$, (ii) it is partially submerged (and filled) with water on the bottom side, and (iii) it is at equilibria with the top side under a non-linear Robin-type boundary condition, where $K_w$ represents a mass-transfer coefficient between the internal saturation at the boundary $\Theta(L, t)$ and the external saturation, that according to \eqref{eq:sat_def} reads $\Theta_\text{ext} = \frac{\theta_\text{ext}-\theta_r}{\theta_s-\theta_r}$.

\begin{remark}
    Here, it is important to notice that the parameter $K_w$ is not rescaled with $c$, but rather intrinsically absorb the role of the diffusion function $D(\Theta)$; as a consequence, the relative importance of boundary exchange versus internal diffusion depends on the calibrated diffusion strength only.
\end{remark}


\subsubsection*{Spatial and temporal discretisation}

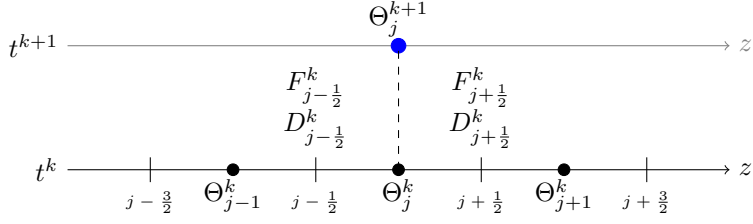
\begin{figure}[t]
    \centering
    \input{figs/gridLayout}
    \caption{
    Finite-difference stencil for the nonlinear diffusion equation.
    At time level $t^k$, nodal values $\Theta_j^k$ are stored at grid points, while interface diffusivities $D_{j\pm\frac12}^k$ and numerical fluxes $F_{j\pm\frac12}^k$ are evaluated at cell interfaces.
    The update of $\Theta_j^{k+1}$ results from the symmetric flux difference across the control volume surrounding node $j$.
    }
    \label{fig:grid}
\end{figure}

The spatial domain $z \in [0,L]$ is discretised into $N_z$ uniform intervals of size $\Delta z = \frac{L}{N_z}$, while time is discretised with step size $\Delta t$, subject to a stability constraint associated with the explicit scheme, \ie the time step is required to satisfy the CFL-type condition
\begin{equation}\label{CFL}
    \frac{\Delta t}{\Delta z^2} \leq \frac{1}{2c\, \max_{[0,1]}{\widehat{D}(\Theta)}} \; ,
\end{equation}
hence ensuring stability of the explicit scheme.
\begin{remark}
    It is important to note that we explicated the role of $c$ within \eqref{CFL} since this highlights that the maximal diffusivity scales linearly with $c$ and, hence, the stability constraint becomes increasingly restrictive for large values of $c$.
\end{remark}

We employ a finite-difference discretisation of the diffusion operator with diffusivity evaluated at cell interfaces.
Let $\Theta_j^k$ denote the numerical approximation of $\Theta(z_j,t_k)$, see also Figure~\ref{fig:grid}.
The explicit update reads
\begin{equation}
    \begin{cases}
    	\Theta_j^{k+1}
    	=
    	\Theta_j^k
    	+
    	\frac{\Delta t}{\Delta z^2}
    	\Bigl[
    	D_{j+\frac12}^k \left(\Theta_{j+1}^k - \Theta_j^k\right)
    	-
    	D_{j-\frac12}^k \left(\Theta_j^k - \Theta_{j-1}^k\right)
    	\Bigr],\\
        \Theta_0^{k+1} = 1\\
        \Theta_{N_z}^{k+1} = \frac{4 \Theta^{k}_{N_z-1} - \Theta^k_{N_z-2} + 2 K_w \Delta z \Theta_{ext}}{3+2K_w\cdot\Delta z},
    \end{cases}
\end{equation}
where the interface diffusivities $D_{j\pm\frac12}^k$ are computed from the neighbouring nodal values simply as
\begin{equation}
    D_{j\pm\frac12}^k = \frac{D_{j}^k + D_{j\pm1}^k}2
\end{equation}
and the two values
\begin{equation}
    F_{j+\frac12} = D_{j+\frac12}^k \left(\Theta_{j+1}^k - \Theta_j^k\right)
    \,,\qquad
    F_{j-\frac12} = D_{j-\frac12}^k \left(\Theta_j^k - \Theta_{j-1}^k\right)
\end{equation}
do represent the numerical fluxes solving \eqref{eq:diff-eq} under a simmetric three-point stencil.
Note that the boundary condition at the top of the specimen reported in the last equation of \eqref{IBC} is discretized using a second order approximation.


\subsubsection*{Computation of the imbibition curve}

For a given parameter vector $\mathbf{p}$, the numerical solution $\Theta(z,t_k)$ defines a model prediction of cumulative water uptake. The simulated imbibition curve is computed as
\begin{equation}
	\widehat{Q}_k^{\mathbf{p}} \equiv \widehat{Q}(t_k; \mathbf{p})
	=
	\rho \int_0^L \Theta(z,t_k)\, dz,
\end{equation}
where $\rho$ is the liquid density (or an appropriate scaling constant depending on normalisation, $\rho = 1\,\sfrac{\text{g}}{\text{cm}^3}$ for water).

In practice, the integral is approximated using the trapezoidal rule:
\begin{equation}
	\widehat{Q}_k
	=
	\rho \sum_{j=0}^{N_z} w_j \, \Theta_j^k \, \Delta z,
\end{equation}
which, under equal weights $w_j$ (but extrema $w_0=\sfrac12=w_{N_{z}}$), simply reads
\begin{equation}
	\widehat{Q}_k
	=
	\rho \frac{\Delta z}{2} \left( \Theta_0^k + 2\sum_{j=1}^{N_z-1} \Theta_j^k + \Theta_{N_z}^k \right)
\end{equation}

This defines the forward operator
\begin{equation}\label{eq:forward}
	\mathbf{p} \;\longmapsto\; \widehat{Q}(t_k; \mathbf{p}) \equiv \widehat{Q}_k^{\mathbf{p}}\,,
\end{equation}
which maps a parameter vector to the corresponding simulated imbibition curve.
The calibration problem then consists of determining the parameter set that minimises the discrepancy between $\widehat{Q}_k^\mathbf{p}$ and the experimental data $Q_k$.


\subsection{Parameter identification problem}

Following \eqref{eq:forward}, the inverse problem consists in determining the parameter vector $\mathbf{p} \in \mathcal{P}$ such that the simulated imbibition curve $\widehat{Q}_k^{\mathbf{p}}$ best reproduces the experimental measurements $Q_k$.
Analogously, this amounts to solving a finite-dimensional optimisation problem over the parameter sets introduced in Section~\ref{sec:parameter_set}. 
For each retention model, the optimisation variables read
\begin{equation}
	\begin{array}{ll}
		\mathbf{p}_{\mathrm{BC}} = (a, b, c, \lambda, \klog),\quad
		&\mathbf{p}_{\mathrm{vG}} = (a, b, c, \lambda, \klog),\\
		\mathbf{p}_{\mathrm{NN}} = (a, b, c, \klog),
		&\mathbf{p}_{\mathrm{BkP}} = (a, b, c, \alpha, \gamma, \klog)
	\end{array}
\end{equation}
where we recall the the definition of $a$ and $b$ from \eqref{eq:a_b} as the relative residual and saturated water content and $c$ being the characteristic diffusion strength of the material, as per \eqref{eq:c}.


\subsubsection*{Cost functional}

To quantify the discrepancy between simulation and experiment, we define the normalised least-squares functional
\begin{equation}
	\mathcal{L}(\mathbf{p})
	=
	\frac{1}{N_T + 1}
	\sum_{k=0}^{N_T}
	\left(
		\frac{\widehat{Q}(t_k; \mathbf{p}) - Q(t_k)}{Q(t_k)}
	\right)^2.
\end{equation}
where the use of relative errors ensures balanced weighting over early and late time regimes of the imbibition curve.

The parameter identification problem can therefore be stated as
\begin{equation}
	\mathbf{p}^\ast
	=
	\arg\min_{\mathbf{p} \in \mathcal{P}}
	\mathcal{L}(\mathbf{p}).
\end{equation}


\subsubsection*{Optimisation strategy}

The minimisation of $\mathcal{L}(\mathbf{p})$ defines a nonlinear and potentially non-convex optimisation problem. 
The strong nonlinearity of the diffusivity $D(\Theta)$ and the structural differences among the retention laws may induce multiple local minima and parameter interactions. 
For this reason, and in order to retain full transparency of the parameter dependencies, we adopt a grid-based search strategy over bounded parameter domains.

For each model, admissible parameter intervals are prescribed based on physical plausibility and ranges reported in the literature. 
A Cartesian product grid is then constructed over these intervals, and the forward simulation is carried out for each parameter combination. 
The loss functional $\mathcal{L}(\mathbf{p})$ is evaluated and the parameter set yielding the minimal value is retained.

Grid search provides a structured and reproducible baseline for parameter identification, allowing direct inspection of parameter sensitivities and interactions. 
However, its computational cost scales exponentially with the number of optimisation variables.
In the present setting, each additional parameter multiplies the total number of forward simulations required.
To mitigate this limitation, we employ an adaptive refinement strategy:
after identifying a candidate optimal region $\mathcal{P}^* \subset \mathcal{P}$ from an initial coarse grid, successive grids are constructed locally within $\mathcal{P}^*$ with progressively finer resolution. 
This procedure concentrates computational effort near promising regions while avoiding exhaustive high-resolution sampling of the entire parameter space.

Despite its computational burden, grid search remains a robust and interpretable reference strategy. 
Alternative optimisation techniques, including Bayesian approaches, heuristic methods, and simulate annealing, have been investigated and systematically compared in our recent works~\cite{Onofri25, Stolfi25}.


\subsubsection*{Implementation details}

All simulations are performed using the explicit finite-difference scheme described above.
The spatial discretisation is fixed to $N_z = 200$ grid points.
The time step $\Delta t$ is selected in accordance with the stability condition \eqref{CFL}, using the maximum diffusivity associated with the current parameter set.

For each candidate parameter vector, the forward problem is solved over the full experimental time horizon and the corresponding imbibition curve is evaluated at the measurement times $t_k$.
The computational effort therefore scales exponentially with the number of parameter explored in the grid search and linearly w.r.t.\ the parameter set sizes.


\section{Numerical results}\label{sec:res}

We now analyse the calibration results obtained for the four retention models introduced in Section~\ref{sec:wrc}. 
The objective is twofold: 
(i) to assess the ability of structurally different diffusivity laws to reproduce the same imbibition dataset, and 
(ii) to investigate how the inferred parameters compare across formulations.

Table~\ref{tab:optim_results} reports the optimal parameter sets identified through progressive grid refinement. 
For each model, the total number of evaluated configurations is indicated, highlighting the computational burden associated with the exploration of multi-dimensional parameter spaces.

\begin{table}[t]
    \centering   
    \input{figs/optimisationTable}
    
    \medskip
    
    \caption{
        Optimization ranges and corresponding best fit for the progressive grid refinement runs (\textbf{It}) on the four water retention models.
        Number of tested configurations is reported in the \textbf{It} column.
    }
    \label{tab:optim_results}
\end{table}


\subsection{Convergence and parameter stability}

Across refinement stages, optimal values remain confined within narrow subregions of the initial search intervals, indicating effective localisation of the minimum. 
The adaptive strategy significantly reduces the need for uniform high-resolution sampling while preserving robustness.

In all models, the relative residual water content $a$ converges to values close to zero, suggesting negligible residual saturation for the tested brick material. 
Conversely, the saturated fraction $b$ systematically approaches values close to unity, reflecting the high effective porosity of the specimen under capillary rise conditions and suggesting a continuous flow in the wet regime (but for BC, where the exponential law requires a larger margin of non-transport in the wet specimen).

\begin{figure}[t]
    \centering
    \includegraphics[width=\linewidth]{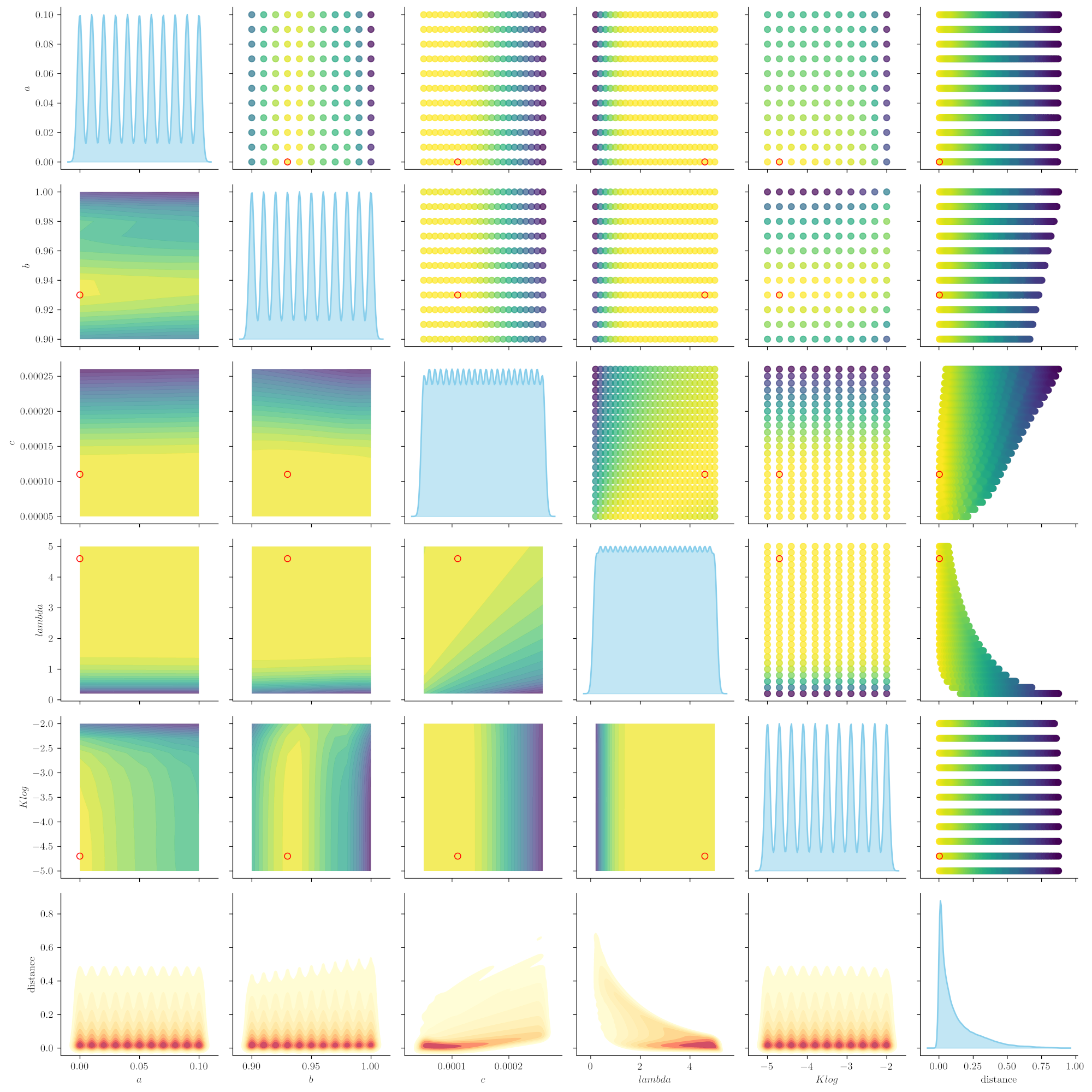}
    \caption{Optimization results for BC model (IT = 2).}
    \label{fig:optim_BC}
\end{figure}

\begin{figure}[t]
    \centering
    \includegraphics[width=\linewidth]{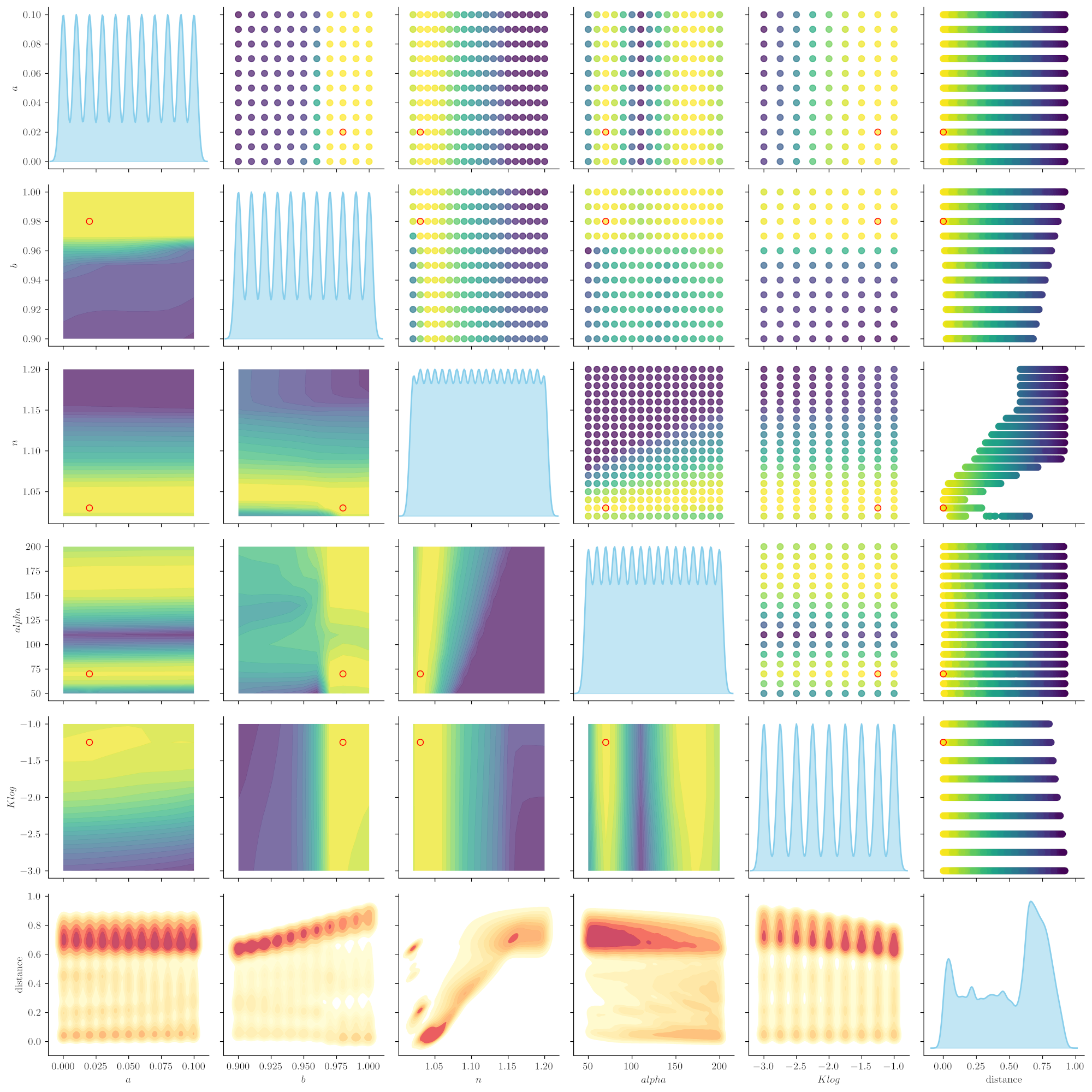}
    \caption{Optimization results for vG model (IT = 1). It is important to note that $\alpha$ is reported instead of $c = \sfrac1\alpha$ and $n$ is reported in place of $\lambda = n-1$.}
    \label{fig:optim_vG}
\end{figure}

\begin{figure}[t]
    \centering
    \includegraphics[width=\linewidth]{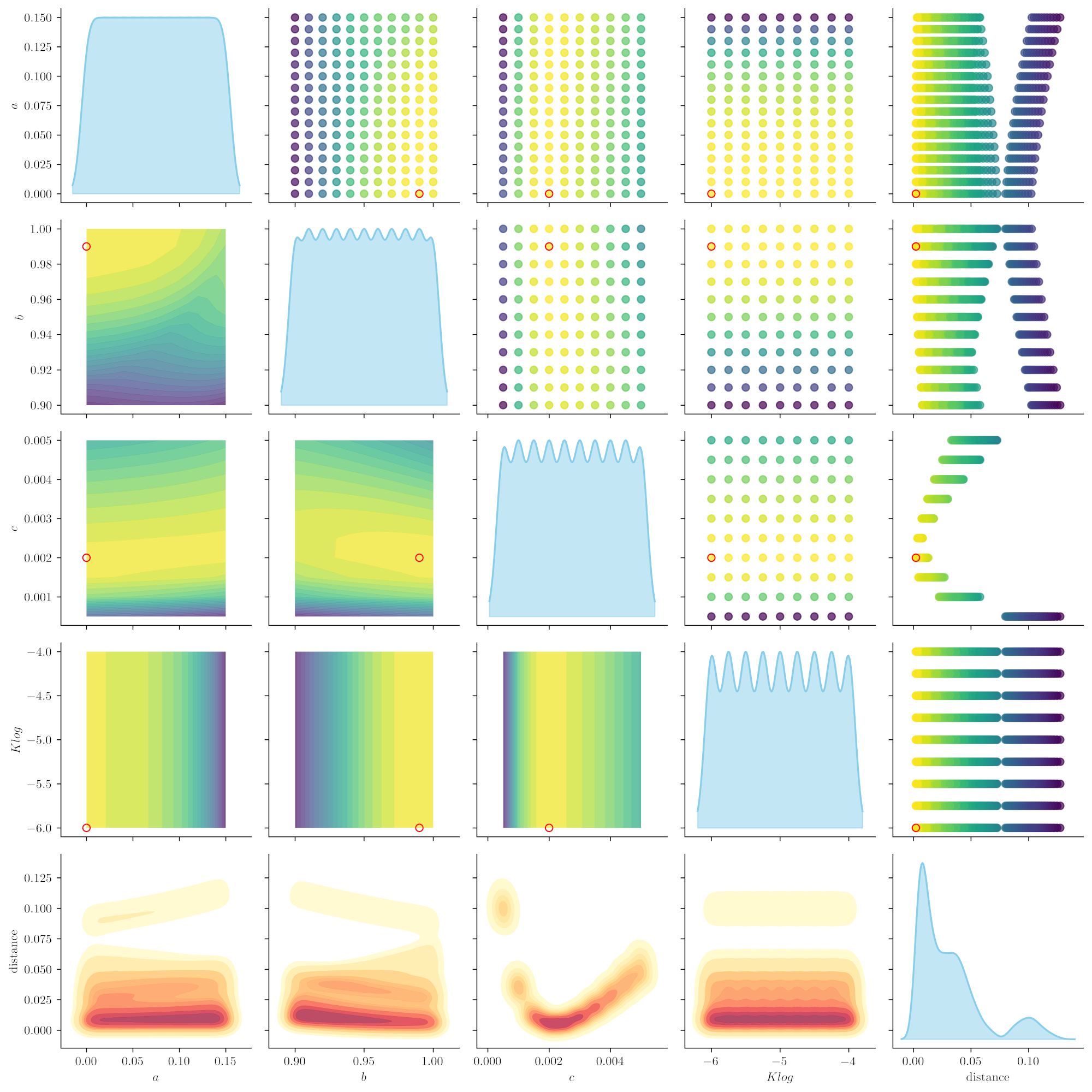}
    \caption{Optimization results for NN model (IT = 1).}
    \label{fig:optim_NN}
\end{figure}

\begin{figure}[t]
    \centering
    \includegraphics[width=\linewidth]{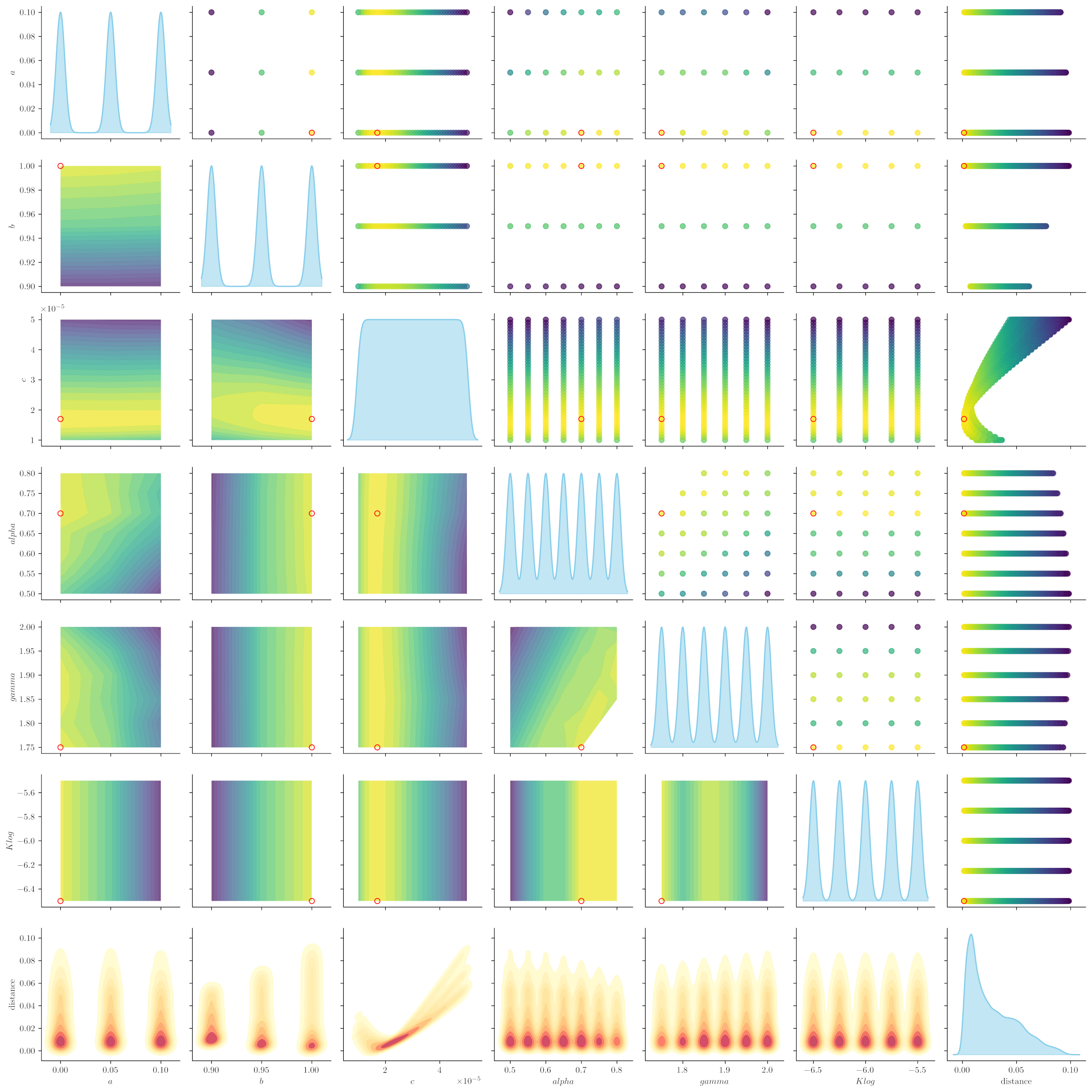}
    \caption{Optimization results for BkP model (IT = 2).}
    \label{fig:optim_BkP}
\end{figure}

The boundary exchange parameter $\klog$ consistently assumes negative values between approximately $-6$ and $-1$, confirming that vapour exchange at the open boundary remains limited compared to capillary inflow at the submerged face, particularly for the NN and the BkP models, where the values are comparable and negligible. Conversely, the vG model appears to favour an exchange rate closer to a non-null-flux condition condition ($\klog=-1.25$).

The localisation process is illustrated in Figures~\ref{fig:optim_BC}--\ref{fig:optim_BkP}, 
which report the loss landscapes at the final refinement stage for each model. 
In all cases, the minimum appears well-defined and confined to a restricted region of parameter space, although a few parameters (like BC's $\lambda$ and $\klog$ in general) offer wider optimality areas. 
BC and vG exhibit relatively sharp valleys along the shape parameter direction, 
whereas NN and BkP display broader plateaus, reflecting stronger parameter compensation effects.
Notably, vG is the only model presenting multiple local minima areas, also highlighted by the shift between the two optimal values for $c = \sfrac1\alpha$.
This behaviour suggests increased parameter interaction within the vG formulation, 
where different combinations of $(\alpha,n)$ can produce comparable effective diffusion profiles. 
Such multiplicity reinforces the non-uniqueness of parameter identification from imbibition data alone.


\subsection{Transport scaling and structural differences}

\begin{figure}[t]
    \centering
    \begin{minipage}{.5\linewidth}
        \centering
        \includegraphics[width=.90\linewidth]{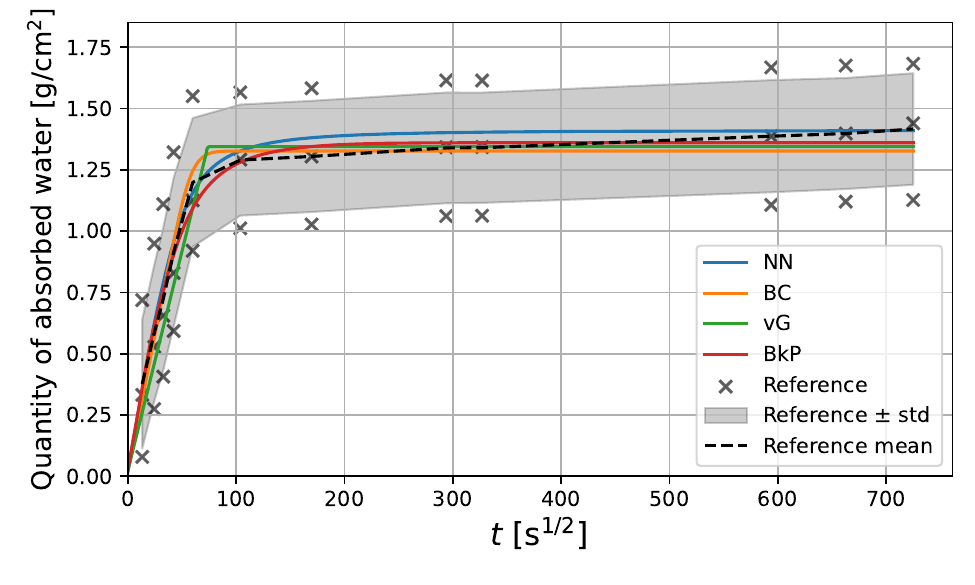}
        \captionof{figure}{Numerical imbibition curves obtained by the four WRC.}
        \label{fig:imb}
    \end{minipage}%
        \begin{minipage}{.5\linewidth}
        \centering
        \includegraphics[width=\linewidth]{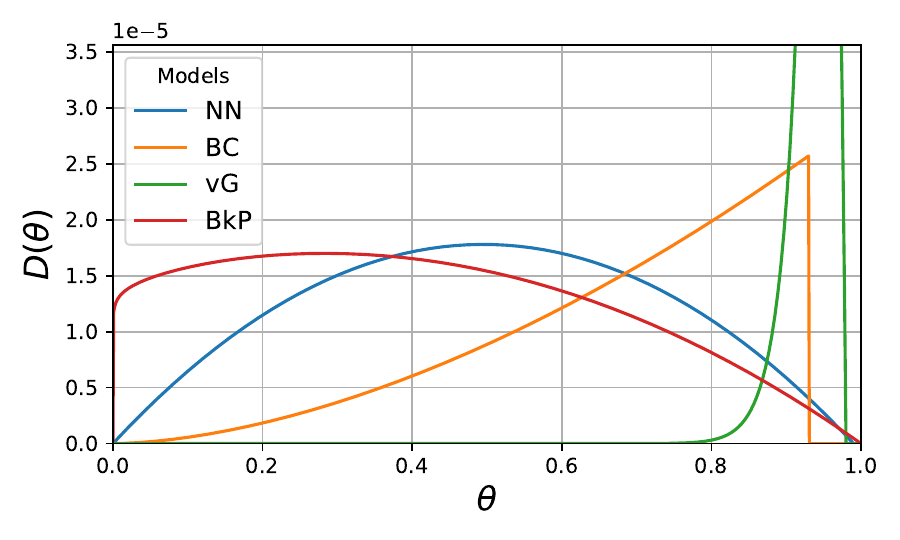}
        \captionof{figure}{Diffusion function of the four WRC. BC and vG are cut to zero for values of $\Theta > b$.}
        \label{fig:D}
    \end{minipage}
\end{figure}

Despite comparable macroscopic fits (Figure~\ref{fig:imb}), the inferred transport scaling parameter $c$ differs markedly across models. 
For the BC model, the optimal value is $c = 1.1\times10^{-4}$, whereas NN and BkP converge to smaller values of order $10^{-5}$. 
For the vG model, the tuned quantity $c = 1.4\times10^{-2}$, is instead much larger, suggesting a different scaling of the diffusivity.

It is important to note that, since $c$ acts as a global scaling factor in the decomposition \eqref{eq:effective_D}, these differences should not be interpreted as direct discrepancies in maximum diffusivity.
Hence, variations in $c$ compensate for intrinsic differences in $\widehat D(\Theta)$.
This effect is clearly visible in Figure~\ref{fig:D}. 

The BC formulation, characterised by a pure exponential law, produces a strongly nonlinear yet monotone diffusion curve. 
The vG model converges to values of $n$ close to unity (i.e.\ small $\lambda=n-1$), indicating behaviour approaching BC-type asymptotics, though with a steeper increase of diffusivity in the intermediate saturation range.

The NN model, constrained by its fixed upside-down parabolic structure, exhibits a symmetric diffusivity profile. 
The reduced scaling $c$ compensates for this structural rigidity. 
Similarly, the BkP model distributes the adjustment between $c$ and the shape parameters $(\alpha,\gamma)$, allowing additional flexibility in both low- and high-saturation regimes.

These results confirm that distinct microscopic assumptions can yield similar macroscopic imbibition curves through compensatory adjustments between transport scaling and functional shape.


\subsection{Macroscopic agreement}

Figure~\ref{fig:imb} shows that all models reproduce the experimental imbibition curve with comparable accuracy. 
Deviations are mainly confined to early times, where steep moisture gradients enhance sensitivity to the diffusion structure and boundary exchange parameter.

The near-indistinguishability of the cumulative uptake curves confirms that imbibition data alone primarily constrain the effective transport intensity, rather than uniquely identifying the detailed constitutive structure of the retention law. 
Different combinations of $(c,\widehat D)$ can therefore produce nearly identical macroscopic responses.


\subsection{Internal moisture dynamics}

\begin{figure}[t]
    \centering
    \includegraphics[width=.47\linewidth]{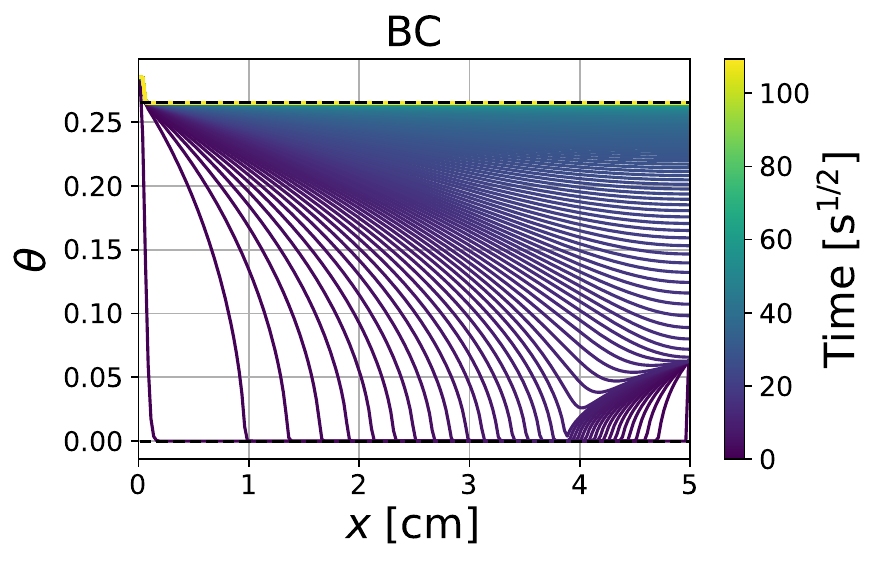}\quad
    \includegraphics[width=.47\linewidth]{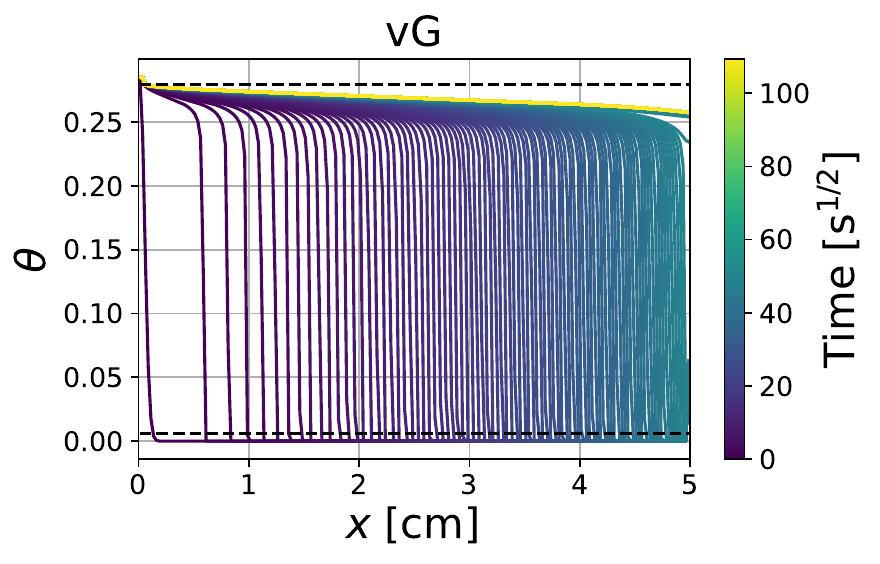}
    
    \includegraphics[width=.47\linewidth]{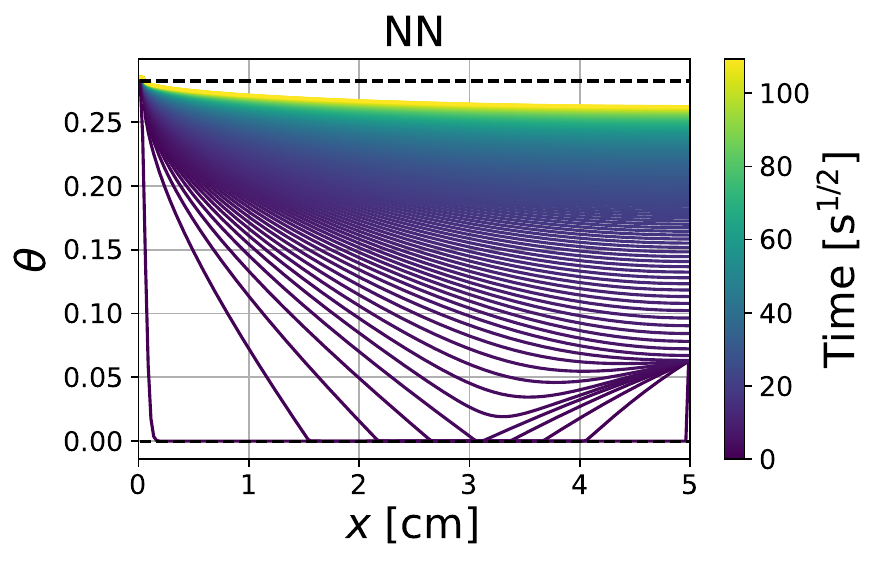}\quad
    \includegraphics[width=.47\linewidth]{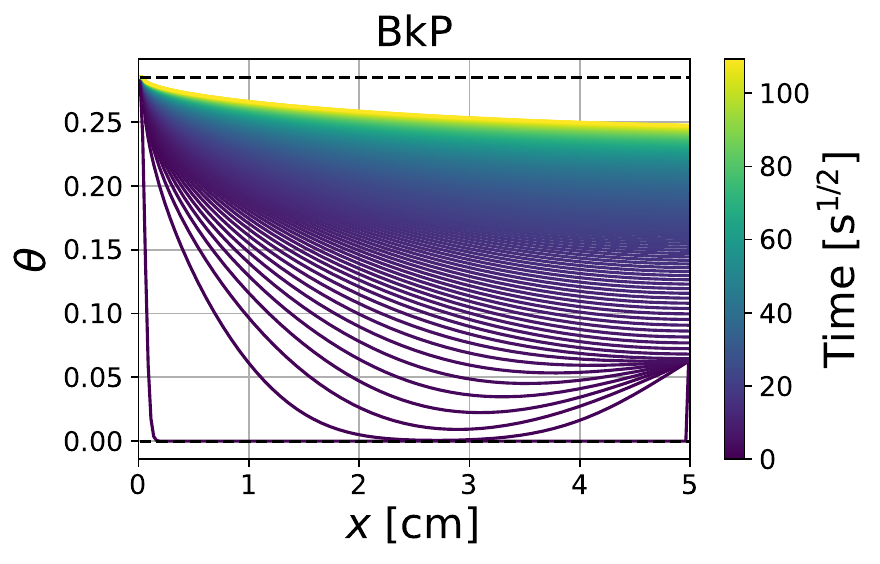}
    \caption{Profile of water content $\theta$ obtained by the different models: BC (top left), vG (top right), NN (bottom left) and BkP (bottom right).}
    \label{fig:theta}      
\end{figure}

Although the macroscopic observable $Q(t)$ is similar across models, internal moisture distributions exhibit noticeable differences (Figure~\ref{fig:theta}). 
The vG retention law produces the sharpest advancing front, due to its rapid increase in diffusivity once a critical saturation threshold is exceeded. 
The BC model yields smoother profiles, since its exponential diffusion law activates progressively over a broader saturation range. 
Conversely, the NN and BkP models generate even smoother transitions, consistent with the regularised parabolic-type structure of their diffusion functions.

\begin{figure}[t]
    \centering\includegraphics[width=.5\linewidth]{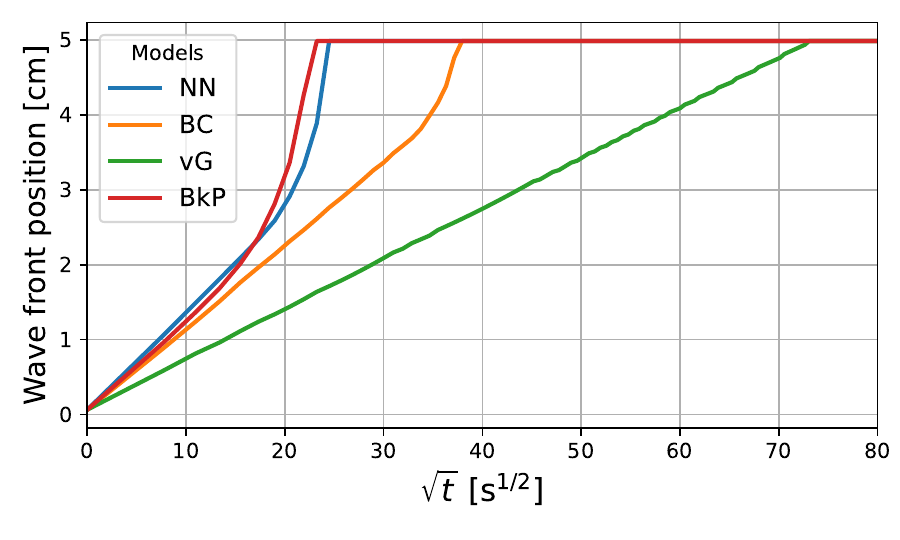}
    \caption{Water fronts obtained numerically with the different models.}
    \label{fig:front}      
\end{figure}

The corresponding front trajectories, reported in Figure~\ref{fig:front}, confirm that similar cumulative uptake may result from distinct propagation mechanisms. In particular, the slowest front advancement is obtained for vG, while BkP and NN show the fastest behavior. 


\subsection{Constitutive comparison}

\begin{figure}[t]
    \centering\includegraphics[width=1.\linewidth]{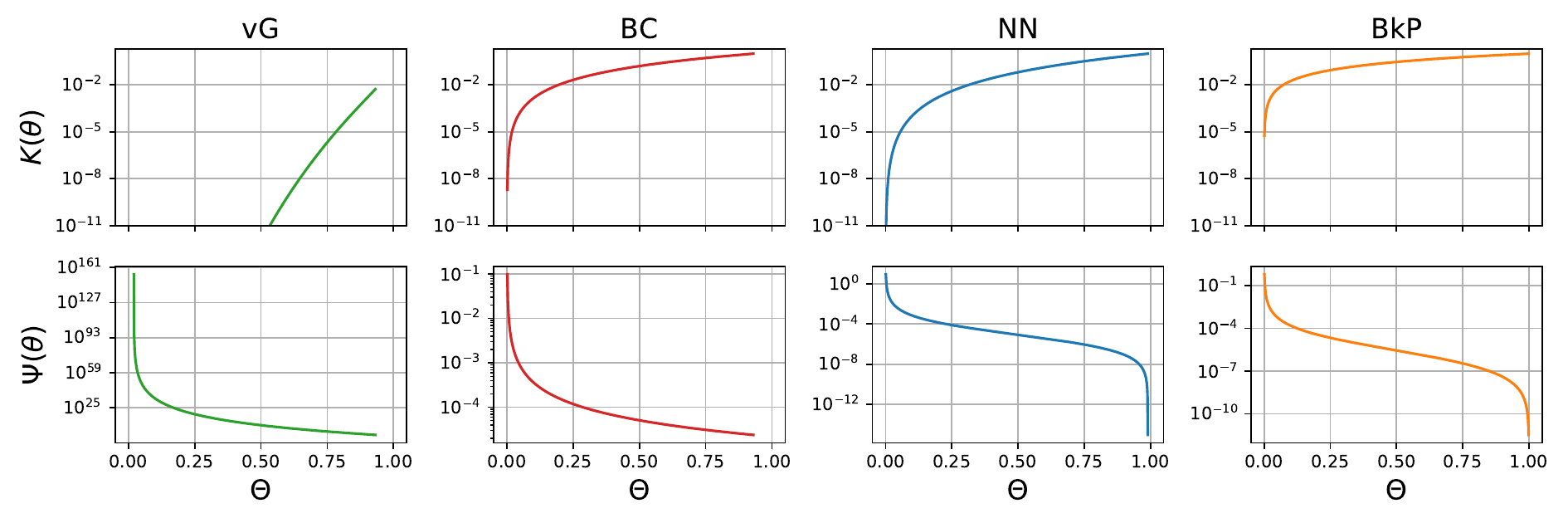}
    \caption{   
        Top Panel: Hydraulic conductivity functions for the different models.
        Bottom Panel: Capillary pressure for the different models.
    }
    \label{fig:fit_all_models}      
\end{figure}

The top and bottom panels of Figure~\ref{fig:fit_all_models} illustrate the calibrated hydraulic conductivity and capillary pressure functions. 
While hydraulic conductivity curves exhibit comparable overall magnitudes --due to the imposed transport scaling-- significant differences emerge in the capillary pressure relations:
BC and vG display hyperbolic-type behaviour with strong curvature at low saturation, 
whereas NN and BkP exhibit smoother transitions, with a two-regime structure reminiscent of hyperbolic-sine behaviour.

This observation highlights an intrinsic identifiability limitation of imbibition experiments: 
the data constrain the effective diffusion coefficient $D(\Theta)$, but do not uniquely determine its decomposition into hydraulic conductivity $K(\Theta)$ and capillary pressure $\psi(\Theta)$. 
Only their combined effect, as reflected in the scaled diffusion function, is structurally observable.


\medskip

Overall, the numerical experiments demonstrate that structurally distinct retention laws can reproduce the same imbibition dataset with comparable accuracy, despite substantial differences in internal moisture dynamics and constitutive behaviour. 
The calibration process primarily constrains an effective transport scaling and the global shape of the diffusion function, while leaving significant flexibility in the decomposition between hydraulic conductivity and capillary pressure. 
These findings motivate a broader discussion on structural identifiability and on the implications of model selection in moisture transport modelling.


\section{Discussion and conclusions}\label{sec:concl}

This study investigated the calibration and comparative behaviour of four nonlinear diffusion models for imbibition experiments within a unified optimisation framework --two physically-based retention laws and two abstract parametrisations. 
The analysis provides insight not only into model performance, but also into structural aspects of parameter identifiability, model flexibility, and computational tractability.

To the best of our knowledge, this is the first systematic comparison of these models calibrated against the same experimental dataset under a consistent numerical and optimisation framework. 
In this sense, the study establishes a first coherent set of parameter configurations for the considered water retention curves that reproduce an identical imbibition experiment, thereby providing a common reference basis for future investigations.

In this regards, the computational framework developed for this study is openly available as part of the author's curated Stoneverse platform\footnote{\url{https://stoneverse.iac.cnr.it}}, see also \cite{Onofri25}. 
By providing open-source access to the models implementation and execution, we aim to facilitate reproducibility, encourage independent validation, and support further methodological developments in the community.


\fakesubsection{Structural identifiability and parameter coupling}
A central outcome of this work concerns the role of the diffusion scaling parameter $c = K_s \psi_b$, which mediates between the diffusion coefficient $D(\Theta)$ and the model-specific shape function $\widehat{D}(\Theta)$ defining the effective diffusion behaviour. 
This factorisation highlights a structural limitation of imbibition experiments: the hydraulic conductivity $K_s$ and the capillary pressure scale $\psi_b$ enter the governing equation multiplicatively and therefore cannot be independently identified from this experimental configuration alone.

Consequently, only their product $c$ is structurally observable. 
This has important implications. 
First, interpreting calibrated values of $K_s$ or $\psi_b$ individually from imbibition data may lead to over-interpretation. 
Second, model comparison must focus on the effective diffusion behaviour rather than on individual hydraulic parameters whose identifiability is intrinsically constrained by the experiment.


\fakesubsection{Model flexibility versus physical interpretability}
The comparison across models reveals a trade-off between physical interpretability and functional flexibility. 
Physically derived retention laws impose structural constraints on the shape of $\widehat{D}(\Theta)$, which enhances interpretability but may limit adaptability to experimental deviations. 
More flexible parametrisations, by contrast, can achieve improved data fitting at the expense of reduced physical transparency.

Differences in predictive behaviour arise primarily from the curvature and degeneracy properties of $\widehat{D}(\Theta)$ near saturation and near dryness, rather than from the overall diffusion scale $c$. 
These regions strongly influence the propagation speed and sharpness of the wetting front. 
Accordingly, a model's ability to reproduce the observed imbibition dynamics depends critically on how it represents the nonlinear structure of the diffusion coefficient, not merely on its magnitude.


\fakesubsection{Implications of the optimisation strategy}
The grid-based optimisation strategy adopted here provides a transparent and reproducible baseline for parameter estimation. 
Its principal advantage lies in its robustness: it does not rely on gradient information and is therefore insensitive to local irregularities in the loss landscape. 
Furthermore, it enables direct inspection of parameter sensitivities and interactions.

However, the computational cost grows exponentially with the number of parameters, reflecting the well-known curse of dimensionality. 
Even with adaptive refinement, the method remains feasible only for moderate parameter dimensions. 
The refinement strategy mitigates this limitation by concentrating resolution in promising regions of the parameter space, but it does not fundamentally alter the scaling behaviour.

Therefore, while grid search serves as a reliable baseline and provides valuable structural insight into the optimisation landscape, it is not intended as a scalable solution for high-dimensional calibration problems. 
Alternative optimisation techniques, including more efficient deterministic and hybrid approaches, have been explored in related studies \cite{Onofri25, Stolfi25}.


\fakesubsection{Sensitivity of imbibition data}
The results indicate that imbibition experiments predominantly constrain the global diffusion scale and the mid-range behaviour of the retention curve, while sensitivity to extreme saturation regimes is comparatively limited. 
This explains why multiple parameter combinations may yield similar loss values when variations mainly affect regions of $\Theta$ that are weakly sampled during the experiment.

An exception is observed for the van Genuchten model, whose steep behaviour near saturation increases sensitivity to parameter variations in that regime. 
More generally, the observed behaviour is consistent with the smoothing character of nonlinear diffusion: the forward operator attenuates high-frequency parameter effects, thereby reducing the effective information content of the data. 
This reinforces the need for complementary experimental configurations if a full hydraulic characterisation is sought.


\fakesubsection{Limitations and outlook}
Several limitations should be acknowledged. 
The analysis is restricted to one-dimensional imbibition under controlled boundary conditions, and model comparison is performed within a deterministic least-squares framework without explicit uncertainty quantification. 
More complex flow scenarios or alternative boundary conditions may alter identifiability properties and parameter sensitivities.

Future work will extend this study along multiple directions. 
First, larger and more diverse datasets will be considered to improve parameter robustness and reduce structural degeneracies. 
Second, probabilistic calibration approaches such as Approximate Bayesian Computation (ABC), as proposed in \cite{Stolfi25}, will be employed to quantify parameter uncertainty and explore posterior distributions. 
Third, Physics-Informed Neural Networks (PINNs), following the approach of \cite{berardi}, will be investigated as a scalable alternative for solving the forward and inverse problems in nonlinear diffusion.


\fakesubsection{Concluding remarks}
Overall, this work clarifies how structural properties of the governing equations influence parameter estimation and model comparison in nonlinear moisture transport. 
Under imbibition conditions, the dominant identifiable quantity is the diffusion scaling parameter $c = K_s \psi_b$, whereas the discriminating power between models lies primarily in the shape of the dimensionless diffusion function $\widehat{D}(\Theta)$.

Beyond the specific models considered, the study highlights the importance of combining structural analysis, careful optimisation strategies, and transparent computational tools when comparing constitutive laws. 
These findings provide a rigorous basis for interpreting calibration results and for guiding future experimental and inversion strategies in porous media flow.


\section*{Authors' contribution}

\noindent
Data Curation: EO.
Conceptualisation: GB, MC.
Formal analysis: EO.
Funding acquisition: GB.
Investigation: GB, EO.
Methodology: GB, MC, EO.
Project administration: GB.
Resources: GB.
Software: EO, MP.
Validation: EO.
Visualisation: EO.
Writing - Original Draft: GB, MC, EO.
Writing - Review \& Editing: GB, EO, MP.
\\
All authors have read and agreed to the published version of the manuscript.
\\
As common practice in mathematics, the authors have not followed the SDC approach; instead, they have sorted the authors alphabetically.


\section*{Acknowledgements}

\noindent Maurizio Ceseri passed away before this manuscript was finalized. We have attempted to present the results of our collaboration in accordance with his high standards. This paper is dedicated to him.

\bigskip

\noindent
G.B.\ and E.O.\ are members of the Gruppo Nazionale Calcolo Scientifico-Istituto Nazionale di Alta Matematica (GNCS-INdAM).\\
G.B.\ is involved in PRIN project MATHPROCULT Prot.\ P20228HZWR, CUP B53D23015940001.\\
G.B., E.O., and M.P.\ are in the PNRR Project H2IOSC CUP B63C22000730005, financed by European Union -- NextGenerationEU PNRR Mission 4, ``Instruction and Research'' -- Component 2 -- Investment line 3.1.\\
G.B.\ is also involved in Project PE0000020 CHANGES -- CUP: B53C22003890006,  NRP Mission 4 Component 2 Investment 1.3, Funded by the European Union - NextGenerationEU under the Italian Ministry of University and Research (MUR).


\bibliographystyle{elsarticle-num}
\bibliography{biblio_arxiv}

\end{document}

%% file: figs/experimentalSetting.tex
\definecolor{botcolor}{HTML}{643A12}
\definecolor{topcolor}{HTML}{CC6E21}

\begin{tikzpicture}
    \draw[fill=cyanM, draw=blueM] (0, 0) rectangle ++(5, .75);

    \shade[bottom color = botcolor, top color = topcolor] (1.0, 0) rectangle ++(3, 3);
    \draw (1, 0) rectangle ++(3, 3);

    \foreach \x in {1.5, 2.0, 2.5, 3.0, 3.5} {
        \draw[-latex, color=cyanM, line width=1.2mm] (\x, 0) -- (\x, 1.5);
        \draw[latex-latex, color=cyanM] (\x, 2.7) -- (\x, 3.3);
    }

    \node at (1.5, 4.0) {$\theta_{ext}$};
    \node at (2.5, 2.0) {$\theta(z, t)$};

    \draw[<->] (0.5, 0.025) -- node[fill=white] {$L$} (0.5, 3);
    \draw[<->] (4.5, 0.025) -- node[fill=cyanM, inner sep=2pt] {\tiny$h_0\sim 0$} (4.5, .75);
    
    \draw[ultra thick] (0, 5) -- (0,0) -- (5, 0) -- (5, 5);
\end{tikzpicture}

%% file: figs/gridLayout.tex
\begin{tikzpicture}[scale=1.1]
    
    \draw[->] (0,1) -- (8,1) node[right] {$z$};
    \draw[->, color=gray] (0,2.5) -- (8,2.5) node[right] {$z$};
    
    \draw (0,2.5) node[left] {$t^{k+1}$};
    \draw (0,1) node[left] {$t^{k}$};
    
    \foreach \x/\lab in {2/j-1,4/j,6/j+1} {
        \filldraw[black] (\x,1) circle (2pt);
        \node[below] at (\x,1) {$\Theta_{\lab}^k$};
    }
    
    \foreach \x/\lab in {1/j-\frac32,3/j-\frac12,5/j+\frac12,7/j+\frac32} {
        \draw (\x,0.85) -- (\x,1.15);
        \node[below] at (\x,0.85) {\tiny$\lab$};
    }
    \foreach \x/\lab in {3/j-\frac12,5/j+\frac12} {
        \node[above] at (\x,1.15) {$D_{\lab}^k$};
        \node at (\x,2.0) {$F_{\lab}^k$};
    }
    
    \filldraw[blue] (4,2.5) circle (2.5pt);
    \node[above] at (4,2.5) {$\Theta_j^{k+1}$};
    
    \draw[dashed] (4,1) -- (4,2.5);
    
\end{tikzpicture}

%% file: figs/optimisationTable.tex
\begin{tabular}{|c|c|c|c|c|c|c|}
    \hline
    \multicolumn{7}{|c|}{\textbf{BC model}}  \\
    \hline  \hline
    \textbf{It} & $a$ & $b$ & $c$ & $\lambda$ & & $\klog$\\
    \cline{1-5}\cline{7-7}
    0 & 0.0 & 0.96 & $2 \times 10^{-4}$ & 10 & & -3.0 \\
    (263538) & $[0, 0.1]$ & $[0.9, 1]$ & $[5\times10^{-5}, 2.5\times 10^{-4}]$ & $[5, 10]$ & & $[-5, -3]$ \\
    \cline{1-5}\cline{7-7}
    1 & 0.0 & 0.93 & $2.6 \times 10^{-4}$ & 11.5 & & -4.0 \\
    $(503118)$ & $[0, 0.1]$ & $[0.9, 1]$ & $[5\times10^{-5}, 2.5\times 10^{-4}]$ & $[5, 15]$ & & $[-4, -2]$ \\
    \cline{1-5}\cline{7-7}
    2 & 0.0 & 0.93 & $1.1 \times 10^{-4}$ & 4.6 & & -4.7 \\
    $(732050)$ & $[0, 0.1]$ & $[0.9, 1]$ & $[5\times10^{-5}, 2.5\times 10^{-4}]$ & $[0.2, 5]$ & & $[-5, -2]$ \\
    \hline
    \multicolumn{1}{c}{}\\
    \hline
    \multicolumn{7}{|c|}{\textbf{vG model}}  \\
    \hline  \hline
    \textbf{It} & $a$ & $b$ & $c^\dagger$ & $\lambda$ (= $n$-1) & & $\klog$\\
    \cline{1-5}\cline{7-7}
    0 & 0.0 & 0.95 & $6.7\times10^{-3}$ & 0.05 & & -2 \\
    (159720) & $[0, 0.1]$ & $[0.9, 1]$ & $[6.7\times10^{-3}, 10^{-1}]$ & $[0.05, 0.40]$ & & $[-5, -2]$ \\
    \cline{1-5}\cline{7-7}
    1 & 0.02 & 0.98 & $1.4\times10^{-2}$ & 0.03 & & -1.25 \\
    (331056) & $[0, 0.1]$ & $[0.9, 1]$ & $[5\times10^{-3}, 2\times10^{-2}]$ & $[0.02, 0.20]$ & & $[-3, -1]$ \\
    \hline
    \multicolumn{1}{c}{}\\
    \hline
    \multicolumn{7}{|c|}{\textbf{NN model}}  \\
    \hline  \hline
    \textbf{It} & $a$ & $b$ & $c^\ddagger$ & \multicolumn{2}{c|}{} & $\klog$\\
    \cline{1-4}\cline{7-7}
    0 & 0.0 & 0.99 & $1.78\times10^{-5}$ & \multicolumn{2}{c|}{} & -5 \\
    (114400) & $[0, 0.15]$ & $[0.9, 1]$ & $[8.9\times10^{-6}, 4.5\times10^{-4}]$ & \multicolumn{2}{c|}{} & $[-5, -2]$ \\
    \cline{1-4}\cline{7-7}
    1 & 0.0 & 0.99 & $1.78\times10^{-5}$ & \multicolumn{2}{c|}{} & -6.0 \\
    $(15840)$ & $[0, 0.15]$ & $[0.9, 1]$ & $[4.5\times10^{-6}, 4.5\times10^{-5}]$ & \multicolumn{2}{c|}{} & $[-6, -4]$ \\
    \hline
    \multicolumn{1}{c}{}\\
    \hline
    \multicolumn{7}{|c|}{\textbf{BkP model}}  \\
    \hline  \hline
    \textbf{It} & $a$ & $b$ & $c$ & $\alpha$ & $\gamma$ & $\klog$\\
    \hline
    0 & 0.0 & 1.0 & $5 \times 10^{-6}$ & 0.1 & 2.0 & -6.0 \\
    (972972) & $[0, 0.1]$ & $[0.95, 1]$ & $[2\times10^{-7}, 5\times 10^{-6}]$ & $[0, 0.1]$ & $[2, 4]$ & $[-6, -4]$ \\
    \hline
    1 & 0.0 & 1.0 & $10^{-5}$ & 0.5 & 1.75 & -6.0 \\
    $(108898)$ & $[0, 0.05]$ & $[0.95, 1]$ & $[10^{-6}, 10^{-5}]$ & $[0, 0.5]$ & $[1.5, 2.5]$ & $[-6, -5]$ \\
    \hline
    2 & 0.0 & 1.0 & $1.7 \times 10^{-5}$ & 0.7 & 1.75 & -6.5 \\
    $(77490)$ & $[0, 0.1]$ & $[0.9, 1]$ & $[10^{-5}, 5\times 10^{-5}]$ & $[0.5, 0.75]$ & $[1.75, 2.0]$ & $[-6.5, -5.5]$ \\
    \hline
\end{tabular}

\smallskip

$^\dagger$ Parameter tuned as $c^{-1}$, corresponding to $\alpha$ from \cite{van-genuchten}.

\smallskip

$^\ddagger$ Parameter tuned as $c\cdot\mu$, as per \cite{Clarelli2010}.

%% file: biblio_arxiv.bib
@article{scarfone,
    title = {A hysteretic water retention model incorporating the soil deformability and its application to rainfall-induced landslides},
    journal = {Computers and Geotechnics},
    volume = {178},
    pages = {106912},
    year = {2025},
    issn = {0266-352X},
    doi = {10.1016/j.compgeo.2024.106912},
    url_ = {https://www.sciencedirect.com/science/article/pii/S0266352X24008516},
    author = {Ehsan Badakhshan and Jean Vaunat and Riccardo Scarfone}
}

@article{berardi,
    title = {Inverse Physics-Informed Neural Networks for transport models in porous materials},
    volume = {435},
    ISSN = {0045-7825},
    url_ = {https://www.sciencedirect.com/science/article/pii/S004578252400882X},
    DOI = {10.1016/j.cma.2024.117628},
    journal = {Computer Methods in Applied Mechanics and Engineering},
    publisher={Elsevier BV},
    author = {Berardi, M. and Difonzo, F. V. and Icardi, M.},
    year = {2025},
    pages = {117628}
}

@article{bracciale,
    title = {Mathematical modelling of experimental data for crystallization inhibitors},
    volume = {48},
    ISSN = {0307-904X},
    url_ = {http://dx.doi.org/10.1016/j.apm.2016.11.026},
    DOI = {10.1016/j.apm.2016.11.026},
    journal = {Applied Mathematical Modelling},
    publisher = {Elsevier BV},
    author = {Bracciale,  M.P. and Bretti,  G. and Broggi,  A. and Ceseri,  M. and Marrocchi,  A. and Natalini,  R. and Russo,  C.},
    year = {2017},
    month = {08},
    pages = {21–38}
}

@article{malte-catania,
  author  = {Bretti, Gabriella and Belfiore, Cristina Maria},
  title   = {Mathematical modelling of water absorption properties for historical lime-based mortars},
  journal = {Int J Geomath},
  year    = {2025},
  volume  = {16},
  number  = {18},
  pages   = {1-18},
  doi = {10.1007/s13137-025-00275-2}
}

@inproceedings{goid,
  title = {Modelling the Effects of Protective Treatments in Porous Materials},
  ISBN = {9783030580773},
  ISSN = {2281-5198},
  _url_ = {http://dx.doi.org/10.1007/978-3-030-58077-3_5},
  DOI = {10.1007/978-3-030-58077-3_5},
  booktitle = {Mathematical Modeling in Cultural Heritage},
  publisher = {Springer International Publishing},
  author = {Bretti,  Gabriella and Filippo,  Barbara De and Natalini,  Roberto and Goidanich,  Sara and Roveri,  Marco and Toniolo,  Lucia},
  year = {2020},
  month = {11},
  pages = {73–83}
}

@article{brooks-corey,
  title = {Properties of Porous Media Affecting Fluid Flow},
  volume = {92},
  ISSN = {2690-3296},
  _url_ = {http://dx.doi.org/10.1061/JRCEA4.0000425},
  DOI = {10.1061/jrcea4.0000425},
  number = {2},
  journal = {Journal of the Irrigation and Drainage Division},
  publisher = {American Society of Civil Engineers (ASCE)},
  author = {Brooks,  R. H. and Corey,  A. T.},
  year = {1966},
  month = {6},
  pages = {61–88}
}

@article{Celia1990,
  title = {A general mass‐conservative numerical solution for the unsaturated flow equation},
  volume = {26},
  ISSN = {1944-7973},
  url_ = {http://dx.doi.org/10.1029/WR026i007p01483},
  DOI = {10.1029/wr026i007p01483},
  number = {7},
  journal = {Water Resources Research},
  publisher = {American Geophysical Union (AGU)},
  author = {Celia,  Michael A. and Bouloutas,  Efthimios T. and Zarba,  Rebecca L.},
  year = {1990},
  month = {7},
  pages = {1483–1496}
}

@inproceedings{Clarelli2010,
  title = {A Mathematical Model for Consolidation of Building Stones},
  _url_ = {http://dx.doi.org/10.1142/9789814280303_0021},
  DOI = {10.1142/9789814280303_0021},
  booktitle = {Applied and Industrial Mathematics in Italy III},
  publisher = {WORLD SCIENTIFIC},
  author = {Clarelli,  Fabrizio and Natalini,  Roberto and Nitsch,  Carlo and Santarelli,  Maria Laura},
  year = {2009},
  month = {9},
  pages = {232–243}
}

@article{fredlund,
    author = {Fredlund, D.G. and Xing, Anqing},
    title = {Equations for the soil-water characteristic curve},
    journal = {Canadian Geotechnical Journal},
    volume = {31},
    number = {4},
    pages = {521-532},
    year = {1994},
    doi = {10.1139/t94-061},
    url_ = {https://doi.org/10.1139/t94-061}
}

@article{ipp,
    title = {Validity limits for the van Genuchten–Mualem model and implications for parameter estimation and numerical simulation},
    journal = {Advances in Water Resources},
    volume = {29},
    number = {12},
    pages = {1780-1789},
    year = {2006},
    doi = {10.1016/j.advwatres.2005.12.011},
    url_ = {https://www.sciencedirect.com/science/article/pii/S0309170805003015},
    author = {O. Ippisch and H.-J. Vogel and P. Bastian}
}

@article{mualem,
  author  = {Mualem, Yechezkel},
  title   = {A new model for predicting the hydraulic conductivity of unsaturated porous media},
  journal = {Water resources research},
  year    = {1976},
  volume  = {12},
  doi = {10.1029/WR012i003p00513},
  issue  = {3},
  pages   = {513-522}
}

@article{Onofri25,
    doi = {10.2312/DH.20253236},
    url_ = {https://diglib.eg.org/handle/10.2312/dh20253236},
    author = {Onofri,  Elia and Bizzarro,  Sofia and Tassa,  Sandro and Czech,  Michela and Bretti,  Gabriella},
    title = {StoneVerse: Models and Methods in Cultural Heritage. The Open-Science Platform for Reproducible Modelling of Stone Decay},
    journal = {Digital Heritage},
    publisher = {The Eurographics Association},
    year = {2025},
    copyright = {Creative Commons Attribution 4.0 International}
}

@article{scelsi,
    title = {Modelling the behaviour of unsaturated non-active clays in saline environment},
    journal = {Engineering Geology},
    volume = {295},
    pages = {106441},
    year = {2021},
    issn = {0013-7952},
    doi = {10.1016/j.enggeo.2021.106441},
    url_ = {https://www.sciencedirect.com/science/article/pii/S001379522100452X},
    author = {Giulia Scelsi and Ayman A. Abed and Gabriele {Della Vecchia} and Guido Musso and Wojciech T. Sołowski},
}

@article{Stolfi25,
  doi = {10.2312/DH.20253258},
  url_ = {https://diglib.eg.org/handle/10.2312/dh20253258},
  author = {Stolfi,  Paola and Onofri,  Elia and Bretti,  Gabriella},
  title = {Estimating Cultural Heritage Processes Using Approximate Bayesian Computation},
  journal = {Digital Heritage},
  publisher = {The Eurographics Association},
  year = {2025},
  copyright = {Creative Commons Attribution 4.0 International}
}

@article{suh,
    author = {Suh, Hyoung Suk and Song, Jun Young and Kim, Yejin and Yu, Xiong and Choo, Jinhyun},
    title = {Data-driven discovery of interpretable water retention models for deformable porous media},
    journal = {Acta Geotechnica},
    volume = {19},
    pages = {3821–3835},
    doi = {https://doi.org/10.1029/2019GL084532},
    url_ = {https://agupubs.onlinelibrary.wiley.com/doi/abs/10.1029/2019GL084532},
    eprint_ = {https://agupubs.onlinelibrary.wiley.com/doi/pdf/10.1029/2019GL084532},
    year = {2024}
}

@article{van-genuchten,
    title = {A Closed‐form Equation for Predicting the Hydraulic Conductivity of Unsaturated Soils},
    volume = {44},
    ISSN = {1435-0661},
    _url_ = {http://dx.doi.org/10.2136/sssaj1980.03615995004400050002x},
    DOI = {10.2136/sssaj1980.03615995004400050002x},
    number = {5},
    journal = {Soil Science Society of America Journal},
    publisher = {Wiley},
    author = {van Genuchten, M. Th.},
    year = {1980},
    month = sep,
    pages = {892–898}
}

@article{vogel,
    title = {Effect of the shape of the soil hydraulic functions near saturation on variably-saturated flow predictions},
    journal = {Advances in Water Resources},
    volume = {24},
    number = {2},
    pages = {133-144},
    year = {2000},
    doi = {10.1016/S0309-1708(00)00037-3},
    url_ = {https://www.sciencedirect.com/science/article/pii/S0309170800000373},
    author = {T. Vogel and M.Th. {van Genuchten} and M. Cislerova}
}

@article{zahasky,
    author = {Zahasky, Christopher and Benson, Sally M.},
    title = {Spatial and Temporal Quantification of Spontaneous Imbibition},
    journal = {Geophysical Research Letters},
    volume = {46},
    number = {21},
    pages = {11972-11982},
    doi = {10.1029/2019GL084532},
    url_ = {https://agupubs.onlinelibrary.wiley.com/doi/abs/10.1029/2019GL084532},
    eprint_ = {https://agupubs.onlinelibrary.wiley.com/doi/pdf/10.1029/2019GL084532},
    year = {2019}
}
